\documentclass{amsart}

\usepackage{mathptmx}

\usepackage{amssymb}
\usepackage{amsmath}
\usepackage{amsfonts}
\usepackage{mathrsfs}
\usepackage{fix-cm}
\usepackage{graphicx}
\usepackage{amscd}
\usepackage{enumerate}
\usepackage{mathptmx}
\usepackage{multirow,booktabs}
\usepackage{hyperref}
\usepackage{tikz}
\usetikzlibrary{arrows,decorations,shapes,shadows,cd}
\usetikzlibrary{babel}
\usetikzlibrary{matrix}

\usepackage{soul}
\sodef\spred{}{.2em}{.9em plus.4em}{1em plus.1em minus.1em}

\usepackage[latin1]{inputenc}
\usepackage[T1]{fontenc}

\usepackage[english]{babel}

\usepackage{bbm}

\usepackage[all]{xy}

\newbox\mybox
\def\overtag#1#2#3{\setbox\mybox\hbox{$#1$}\hbox to
  0pt{\vbox to 0pt{\vglue-#3\vglue-\ht\mybox\hbox to \wd\mybox
      {\hss$\ss#2$\hss}\vss}\hss}\box\mybox}
\def\undertag#1#2#3{\setbox\mybox\hbox{$#1$}\hbox to 0pt{\vbox to
    0pt{\vglue#3\vglue\ht\mybox\hbox to \wd\mybox
      {\hss$\ss#2$\hss}\vss}\hss}\box\mybox}
\def\lefttag#1#2#3{\hbox to 0pt{\vbox to 0pt{\vss\hbox to
      0pt{\hss$\ss#2$\hskip#3}\vss}}#1}
\def\righttag#1#2#3{\hbox to 0pt{\vbox to 0pt{\vss\hbox to
      0pt{\hskip#3$\ss#2$\hss}\vss}}#1}
\let\ss\scriptstyle

\def\Circ{\hbox{\hss$\circ$\hss}}
\def\Vdots{\hbox{$\vdots$}}
\def\splicediag#1#2{\xymatrix@R=#1pt@C=#2pt@M=0pt@W=0pt@H=0pt}

\newcommand\lineto{\ar@{-}}
\newcommand\dashto{\ar@{--}}
\newcommand\dotto{\ar@{.}}

\newtheorem{theorem}{Theorem}[section]
\newtheorem{lemma}[theorem]{Lemma}
\newtheorem{proposition}[theorem]{Proposition}
\newtheorem{corollary}[theorem]{Corollary}
\newtheorem{algorithm}[theorem]{Algorithm} 

\theoremstyle{definition} 
\newtheorem{definition}[theorem]{Definition} 
\newtheorem{example}[theorem]{Example}

\newcommand{\C}{\mathbbm{C}}

\newcommand{\R}{\mathbbm{R}}
\newcommand{\Z}{\mathbbm{Z}}
\newcommand{\K}{\mathbbm{K}}
\newcommand{\z}{\mathcal{Z}}

\newcommand{\rank}{\operatorname{rank}}

\newcommand{\lcm}{\operatorname{lcm}}

\newcommand{\morf}[4][\to]{ #2 \colon #3 #1 #4}

\renewcommand{\phi}{\varphi}
\renewcommand{\epsilon}{\varepsilon}

\begin{document}

\bibliographystyle{alpha}

\title{Classifying metrically conical rational singularities}
\author{Helge M\o{}ller Pedersen}
\date{}

\maketitle

\begin{abstract} We will in this article classify which rational complex
  surface singularities are metrically conical, that is have the
  bilipschitz type of their inner metric equivalent to that of a
  cone. This is done using the thick-thin decomposition for the inner
  geometry of a normal complex surface singularity and properties of
  rational singularities.

The result is a family of singularities depending only on their
topology which is given as a $3$-parameter family.


\end{abstract}

\section{Introduction}

Bilipschitz geometry is the study of the bilipschitz classes of the
inner and the outer metric. The latter is defined by choosing an
embedding $X\hookrightarrow \R^n$ and then just taking the restriction of the
Euclidean distance function. The inner geometry is defined as the path
metric on $X$, that is
\begin{align*}
d_{in}(u,v):=\inf_{\gamma} \Bigl\{ length_{\R^n}(\gamma)\ \Big\vert\
\morf{\gamma}{[0,1]}{X}\text{ where }\gamma(0)=u\text{ and }\gamma(1)=v\Bigr\}.
\end{align*}A priori the bilipschitz class of these metrics depends on
the embedding, but if $(X,0)$ is the germ of an algebraic or analytic
variety then it turns out that the bilipschitz class only depends on
the analytic type, and not on the chosen embedding.

It is clear that if $(X,0)$ and $(Y,0)$ are bilipschitz equivalent then
they are homeomorphic, and it is also not hard to see that if $(X,0)$
and $(Y,0)$ are bilipschitz equivalent with respect to the outer
metric, then they are also bilipschitz equivalent with respect to the
inner metric. This implies that we can see the bilipschitz
classification of singularities as sitting between the analytic
classification and the topological classification. This has
advantages since due to a result of Mostowski
\cite{MostowskiLsatrtifications1} and Parusi\'nski 
\cite{ParusinskiLstratifications2} the
bilipschitz classification does not have continuous moduli in
contrast to the analytic classification.

Any singular point $(X,0)$ of a real or complex variety satisfies the
Conical Structure Theorem. That is given an embedding
$(X,0)\hookrightarrow (\K,0)$ (here $\K=\R$ or $\C$) there exists an
$\epsilon_0$ such that the homeomorphism type of $X\cap S^N_\epsilon$
for $0<\epsilon<\epsilon$ does not depend on any choices made, here
$S^N_\epsilon$ is the sphere of radius $\epsilon$. The homeomorphism
type of this intersection is called the link of $(X,0)$.

So the first question that comes to mind is, are the bilipschitz
classes of $(X,0)$ just cones over geometries on the link? The answer
to this is in general no. It already follows from the work of Pham and
Teissier \cite{phamteissier} that the outer classification of complex
plane curves is the same as the embedded topological
classification. But for complex curves the answer for the inner metric
is yes. The germ of an irreducible complex curve is inner bilipschitz
equivalent to a disc with the standard metric.

The study of the bilipschitz geometry of complex surface singularities
began by Birbrair and Fernandes asking if this was also the case for
surface germes. They found the first counterexample in
\cite{birbrairfernandes}. They were then joined by Neumann and in
\cite{birbrairfernandesneumann1,birbrairfernandesneumann2,birbrairfernandesneumann3}
they found that this was in general not the case. They also found
examples of surface singularities homeomorphic but not inner
bilipschitz equivalent, thereby showing that the inner bilipschitz
classification is stronger than the topological classification. In
light of these findings the following definition makes sense.
\begin{definition}
A singularity $(X,0)$ is called \emph{metrically conical} if its inner
geometry is bilipschitz equivalent to a metric cone over its link.
\end{definition}

In \cite{thickthin} Birbrair, Neumann and Pichon define what they call a
thick-thin decomposition of $(X,0)$ which is invariant under inner bilipschitz
semi-algebraic homeomorphisms. $(X,0)$ is then metrically conical if and only
if it consists of one thick piece and no thin pieces. A consequence of this is
that the link of a metrically conical singularity $(X,0)$ has to be a
Seifert fibred manifold.

Starting from this result we began studying if one can classify which
complex surface singularities are metrically conical, especially if
there are classes of surface singularities where one can give
explicit families of metrically conical singularities. A natural
place to start is the rational singularities. They are identified by
having geometric genus $0$. But, as we will see in Section \ref{preliminaries},
whether a singularity is rational or not is determined by its
topology. The main result of this article will be the following
theorem which indeed gives an explicit description of all rational
singularities that are metrically conical.
\begin{theorem}\label{maintheorem} 
Let $(X,0)$ be a rational surface singularity. Then $(X,0)$ is
metrically conical if and only if it is in one of the two following
families:
\begin{enumerate}[(a)]
\item The link of $(X,0)$ is Seifert fibred with $n>2$ singular
  fibres and Seifert invariants
  $\Bigl(0,\bigl(1,-(n-1)\bigr),\bigl(p,p-1\bigr) ,\dots,
  \bigl(p,p-1\bigr), \bigl(p,p-q\bigr), \bigl(p,q\bigr)\Bigr)$ where
  $p$ and $q$ are positive integers with $p> q$ and
  $\gcd(p,q)=1$.\label{threeormoreseifertfibrescase}

\item The exceptional divisor of the minimal good resolution of
  $(X,0)$ consists of a single smooth rational curve. Equivalently, the
  link is a circle bundle over $S^2$.\label{singlecurvecase}
\end{enumerate}
\end{theorem}

As we will see later, for rational singularities the inner
classification coincides with the topological classification, so it is
no surprise this classification relies on the topology. In fact as we
will see in Section \ref{preliminaries} what we really prove is the
following equivalent theorem.

\begin{theorem}\label{nonmetricversion}
Let $(X,0)$ be a normal complex surface singularitiy with Siefert
fibered link. Then $(X,0)$ is
rational and the strict transform of a generic hyperplace section only
intersects the central node of the minimal good resolution of $(X,0)$
if and only if either \eqref{threeormoreseifertfibrescase} or
\eqref{singlecurvecase} holds.
\end{theorem}

If we are in case \eqref{singlecurvecase} of Theorem \ref{maintheorem} then the
simplest example is when the unique curve in the exceptional divisor
has self intersection number $-2$. But this is the $A_1$ simple
singularity which can be given by the equation $x^2+y^2+z^2=0$. In
case \eqref{threeormoreseifertfibrescase} if $n=3$, $p=2$ and $q=1$ then the Seifert invariants are
$\bigl(0,(1,-2),(2,1),(2,1),(2,1)\bigr)$ which is the simple $D_4$
singularity. It is not hard to see that non of the other simple $A_n$,
$D_n$, $E_6$, $E_7$ and $E_8$ singularities are in the families
described in Theorem \ref{maintheorem}. Since these are all the rational
singularities that are hypersurfaces, it follows that only $A_1$ and
$D_4$ are metrically conical rational hypersurfaces. This fact was
already observed in \cite{thickthin}.

Since the classification of the inner geometries on normal complex
surfaces was done in \cite{thickthin}, there has been a lot of work on
bilipschitz
geometry. We will not mention everything, but of note is the proof that
outer bilipschitz regularity implies smoothness for complex
singularities by Birbrair, Fernandes, L\^e and Sampaio
\cite{lipschitzregularity}, work on showing that Lipschitz triviality
is the same as Zariski equisingularity by Neumann and Pichon
\cite{zariski} and Parusi\'nski and P\u aunescu
\cite{parusinskipaunescu}, the introduction of \emph{moderately
  discontinuous homology} by Bobadilla, Heinze, Pe Pereira and Sampaio
\cite{BobadillaHeinzePereiraSampaio}, proof that for hypersurfaces in
$(\C^3,0)$ and $(\C^4,0)$ the multiplicity is a bilipschitz
invariant for the outer metric by Fernandes and Sampaio
\cite{fernandessampaioMultdim23}  and
Sampaio \cite{sampaioMatricZarCOndim4},
and a lot of work on the notion
of being \emph{Lipschitz normally embedded}
\cite{birbrairMostowski,BirbrairMendes,LNEimpliesLNEtangentcones,lneinthespaceofmatrices,charecterizationofLNE,LNEofminimalsingualrities,MisevPichon}. 

The question of being metrically conical was on the other hand more or
less dormant
until recently when Kerner and Mendes considered it in their study
on fast cycles
\cite{kernermendesfastvanishingloops,kernermendefastvanishingcycles}.
It was these works that inspired us to look at this question again,
and though we do not actually use any of their results, it did get us
to consider weighted homogeneous rational singularities and thereby
discovering the work of Pinkham \cite{pinkhamcstaraction}, which turned
out to be the missing ingredient in making a proof of Theorem
\ref{maintheorem}.

The article is organized as follows. In section \ref{preliminaries} we
introduce the thick-thin decomposition for the inner geometry of a
normal complex surface, Artin's fundamental cycle and discuss
the properties of rational and weighted homogeneous singularities
needed. In section \ref{suffeciency} we prove that the two families given in
\ref{maintheorem} are indeed both rational and metrically
conical. Finally in section \ref{necessity} we prove that any rational and
metrically conical singularity must belong to one of the two families.

\subsection*{Acknowledgements}
The author would like to thank Anne Pichon who originally proposed
this question to me when I went to work with her as a postdoc in 2011, and
would like to thank the anonymous reviewer that pointed out that the
family I originally found those many years ago was not complete.

\subsection*{A.I. Statement}

The author used generative A.I. non of the other to help find the proof of
Proposition \ref{commonseifertdegree}, the proof is simple, and was
verified and extended by the author to the general setting we use it in. The
author assumes full responsibility for the correctness of everything in
this article.

\section{Preliminaries}\label{preliminaries}

In this section we will introduce the tools we are going to use in
deciding which surface singularities are both rational and metrically
conical. We will introduce different tools in different subsections,
and while the tools in different subsections are independent from one
another, the notation introduced in previous sections will be used
throughout the article.

\subsection{Thick Thin decomposition}

In \cite{thickthin} Birbrair, Neumann and Pichon gave a decomposition
of a normal complex surface germ $(X,0)$ into semi-algebraic pieces
that are either ``thick'' or ``thin''. We will not repeat the precise
definition here, but heuristically a piece is ``thick'' if it is a
``fattened'' up metrically conical semi-algebraic germ (The
``fattening'' is to make gluing work). A ``thin'' piece is on the other
hand a piece whose ``diamater'' collapses faster than linearly preventing
conicalness. This intuition implies the following, which is Corollary
7.11 of \cite{thickthin}

\begin{proposition}\label{thickthinclassyficationofMC}
$(X,0)$ is metrically conical if and only if there are no thin
pieces. Equivalently, $\Gamma$ has only one node, and it is the unique
$\mathcal{L}$-node.
\end{proposition}Now the second version, which is the version we will
actually use, has some words and symbols that might require some
explanation. First $\Gamma$ is a dual graph of a resolution
$\morf{\pi}{(\widetilde{X},E)}{(X,0)}$ that is good and factors
through both the blow-up of the origin and the Nash-Transform of
$(X,0)$. It can be obtained from the minimal good resolution by
blowing up more points. Remember that a resolution is good if the
exceptional divisor is a normal crossing divisor and it is minimal
if any exceptional curve $E$ that corresponds to a vertex of valence $1$ or $2$
has self intersection number $E^2 \leq -2$.

The condition on the resolution secures that it resolves the base
points of both a generic hyperplane section $\morf{h}{(X,0)}{\C,0)}$
and the polar curve $\Pi$ of $(X,0)$. This implies that the strict
transforms $h^*$ and $\Pi^*$ do not pass through singular points of
the exceptional divisor. More precisely if $E=\bigcup E_i$ then they do
not pass through points in $E_i\cap E_j$. Then any vertex of $\Gamma$
that corresponds to an exceptional curve $E_i$ with $E_i\cap h^*\neq
\emptyset$ is called an $\mathcal{L}$-node. A vertex with $E_i\cap
\Pi^* \neq \emptyset$ is called a $\mathcal{P}$-node. Finally a
vertex with valence more than $2$ or corresponding to a curve with genus
more than $0$ is called a node (it can be a
$\mathcal{L}$-node, $\mathcal{P}$-node or neither.)

This should explain the second line of Proposition
\ref{thickthinclassyficationofMC}. It should also make clear that the
dual resolution graph is star-shaped with the central vertex its unique
$\mathcal{L}$-node. It follows from Proposition 6.1 of
\cite{thickthin} that if $\Gamma$ is star-shaped and the central vertex
is the only $\mathcal{L}$-node, any $\mathcal{P}$-node is also the
central vertex and
hence $(X,0)$ is metrically conical. So we can assume that $\Gamma$
has the following form.
$$\splicediag{8}{30}{
  &&&&\\
 &&& \overtag\Circ{-e_{1k_1}}{10pt}\lineto[dl]\\
  &
  &\overtag\Circ{-e_{1{k_1-1}}}{10pt}\dashto[ddl]&
  \undertag\Circ{-e_{2k_2}}{0pt}\lineto[dl]\\    
  &&\undertag\Circ{-e_{2{k_2-1}}}{0pt}\dashto[ddl]&&\\
 &\overtag\Circ{-e_{11}}{10pt}\lineto[ddl]
  &&&\\
&\undertag\Circ{-e_{21}}{-4pt}\lineto[dl]&&&\\
  \undertag{\overtag\Circ{-a}{10pt}}{[g]}{1pt}\lineto[ddr]
  &\Vdots&&&  \\
&&&\\
  &\overtag\Circ{-e_{n1}}{10pt}\dashto[ddr]
  &&\\
  &&&\\
  &&\overtag\Circ{-e_{n(k_n-1)}}{10pt}\lineto[dr]&&\\
    &&&\overtag\Circ{-e_{nk_n}}{10pt}
       \\
  &&\text{Figure 1.}&&\hbox to 0 pt{~,\hss} }$$
Now having a star-shaped resolution graph like is the same as the
link of $(X,0)$ being Seifert fibred with Seifert invariants
$\bigl(g,(1,-a),(p_1,q_1),\dots,(p_n,q_n)\bigr)$ where $(p_i,q_i)$ are
given as the continued fractions:
\begin{align*}
p_ i/q_i = e_{i1}-\cfrac{1}{e_{i2}-\cfrac{1}{e_{i3}-\dots}}.
\end{align*}Notice that there do not exist a unified notion of what the
Seifert invariants are, so this is the definition we will use. The other
definition has the continued fraction be equal to $p_i/(p_i-q_i)$.

\subsection{Maximal Ideal and Fundamental Cycles}

It follows from the last section that what is needed to prove that
$(X,0)$ is metrically conical, is to detect where the
$\mathcal{L}$-nodes are in the resolution $\Gamma$. That is to
understand the strict transform $h^*$ of a generic hyperplane
section. Now the total transform of the generic linear function $h$
consists of two parts, the strict transform $h^*$ and a cycle
$Z_{max}$ called the \emph{maximal ideal cycle} supported on the
exceptional divisor. Since $Z_{max}$ is
the cycle coming from the pullback of a holomorphic function it follows that
$Z_{max}\cdot E_i \leq 0$ for any exceptional curve $E_i$, and that
the inequality is strict exactly when the strict transform $h^*$
intersects $E_i$.

One can define the set $\z_{top}$ of \emph{topological cycles} as the set of
effective divisors with support on the exceptional divisor, such that
$Z\cdot E_i \leq 0$ for all exceptional curves. The set of topological
cycles is partially ordered if we say that $\sum a_iE_i$ is less than
$\sum b_iE_i$ if $a_i\leq b_i$ for all $i$. Then there is a minimal
topological cycle $Z_{min}$ called the \emph{minimal cycle} or Artin's
\emph{fundamental cycle}. It is clear that the fundamental cycle only
depends on the combinatorics of the resolution graph, and hence the
topology of $(X,0)$.

It is also clear that $Z_{min}\leq Z_{max}$, and
if equality is true then we could use $Z_{min}$ instead of $Z_{max}$
to find the $\mathcal{L}$-nodes. This would be advantageous since we
can find $Z_{min}$ given $\Gamma$ using the
following algorithm due to Laufer \cite{lauferrational}
\begin{algorithm}[Laufer's Algorithm]\label{laufersalgorithm}
Let $E_1\cup\dots\cup E_n$ be the exceptional divisor of a good resolution of
$(X,0)$. Then let $Z_1=E_j$ for some $j$. If there exists an $E_k$ such
that $Z_j\cap E_k > 0$, then define $Z_{j+1}=Z_j+E_k$. This procedure
stops, and the final cycle obtained is $Z_{min}$. 
\end{algorithm}

The first exceptional curve chosen in the algorithm can be arbitrary,
so we can always choose the central curve $E$ if the resolution diagram is
star-shaped. We can then ignore the central curve until there are no
other curves that need to be added. We will denote this cycle as
$Z_1$. Notice that the multiplicity of $E$ in $Z_1$ is 1, and it is
not hard to see that the same is true for all other curves, i.e.\ $Z_1
= E + \sum E_{ij}$. If $Z_1\cdot E \leq 0$ then the algorithm
stops. If not we need to add $E$, we then will add any other curves
needed until it is only the central node that might not satisfy the
inequality and denote this by $Z_2$. We will continue this process and
denote the successive cycles obtained this way as $Z_k$.

Notice that $Z_k$ has the following properties:
\begin{enumerate}
\item The multiplicity of $E$ in $Z_i$ is $i$.
\item The multiplicity of $E_{ij}$ is less than or equal to $i$ for all
  non central curves.
\item $Z_i\cdot E_{ij} \leq 0 $ for all non central curves $E_{ij}$.
\item $Z_{min}=Z_k$ for some $k$.
\end{enumerate}
We will continue to use the notation $Z_k$ for this cycle throughout
the article.

Notice that if the graph is a line and all the genera are $0$, then any
vertex can be chosen as the central one. Since $E^2\leq -2$ we get
$Z_1\cdot E \leq E^2 +2 \leq 0$ for any curve $E$ because there are at
most two neighbouring vertices. This gives:

\begin{corollary}\label{linegraph}
If the minimal good resolution graph is a line, then $Z_{min}=Z_1$.
\end{corollary}

\begin{example}
Consider the dual resolution graph $\Gamma$ below.
\begin{center}
\begin{tikzpicture}
\draw[thin](0,0)--(1,0);
\draw[thin](1,0)--(2,0);
\draw[thin](2,0)--(2,-1);
\draw[thin](2,-1)--(2,-2);
\draw[thin](2,0)--(3,0);
\draw[thin](3,0)--(4,0);
\draw[thin](4,0)--(5,0);
\draw[thin](4,0)--(5,0);
\draw[thin](5,0)--(6,0);

\draw[fill=white](0,0)circle(2pt);
\draw[fill=white](1,0)circle(2pt);
\draw[fill=white](2,-1)circle(2pt);
\draw[fill=white](2,-2)circle(2pt);
\draw[fill=white](2,0)circle(2pt);
\draw[fill=white](3,0)circle(2pt);
\draw[fill=white](4,0)circle(2pt);
\draw[fill=white](5,0)circle(2pt);
\draw[fill=white](6,0)circle(2pt);

\node(a)at(0,0.35){$-3$};
\node(a)at(1.60,-1){$-3$};
\node(a)at(1.60,-2){$-2$};
\node(a)at(1,0.35){$-2$};
\node(a)at(2,0.35){$-2$};
\node(a)at(3,0.35){$-2$};
\node(a)at(4,0.35){$-2$}; 
\node(a)at(5,0.35){$-2$};
\node(a)at(6,0.35){$-2$};

\node(a)at(0,-0.35){$v_{22}$};
\node(a)at(2.4,-1){$v_{31}$};
\node(a)at(2.4,-2){$v_{32}$};
\node(a)at(1,-0.35){$v_{22}$};
\node(a)at(2.3,-0.35){$v$};
\node(a)at(3,-0.35){$v_{11}$};
\node(a)at(4,-0.35){$v_{12}$}; 
\node(a)at(5,-0.35){$v_{13}$};
\node(a)at(6,-0.35){$v_{14}$};

\end{tikzpicture}
\end{center}
The Seifert invariants are
$\bigl(0,(1,-2),(5,4),(5,3),(5,2)\bigr)$. Let $E_k$ correspond to the
vertex $v_k$. Then running Laufer's algorithm
on $\Gamma$ will give:
\begin{align*}
Z_1 &= E +E_{11} +E_{12} +E_{13} +E_{14} +E_{21} +E_{22} +E_{31}
+E_{32} \\
Z_2 &= 2E +2E_{11} +2E_{12} +2E_{13} +E_{14} +2E_{21} +E_{22} +E_{31}
+E_{32} \\
Z_3 &= 3E +3E_{11} +3E_{12} +2E_{13} +E_{14} +2E_{21} +E_{22} +2E_{31}
+E_{32} \\
Z_4 &= 4E +4E_{11} +3E_{12} +2E_{13} +E_{14} +3E_{21} +E_{22} +2E_{31}
+E_{32} \\
Z_5 &= 5E +4E_{11} +3E_{12} +2E_{13} +E_{14} +3E_{21} +E_{22} +2E_{31}
+E_{32}.
\end{align*}Notice that $Z_5= Z_{min}$ and that $Z_5\cdot
E_{ij}=0$ for $j>1$. We also see that $Z_5=Z_4 +E$ and  $Z_4\cdot
E_{ij}=0$. Moreover $Z_i\cdot E = 1$ if $i<5$. These facts turn out
to be general and will be important in our later proofs. 
\end{example}

\subsection{Rational Surfaces Singularities}

Now $Z_{max}$ is not in general equal to $Z_{min}$, but in certain
families that is true. So we will consider one such family.

\begin{definition}
A normal surface singularity is called \emph{rational} if its geometric genus
$p_g:=\rank H_1(\widetilde{X},\mathcal{O}_{\widetilde{X}})$ is $0$.  
\end{definition}

This definition is of course not very useful for us, so we will
instead use the following criterion due to Laufer \cite{lauferrational}

\begin{theorem}[Theorem 4.2. \cite{lauferrational}]\label{laufercriterion}
If all the curves in the exceptional divisor have genus $0$, the
resolution graph is a tree and for every step in Laufer's Algorithm
\ref{laufersalgorithm} $Z_j\cdot E_k \leq 1$ for all $E_k$. Then
$(X,0)$ is a rational singularity.
\end{theorem}This shows that being rational is not only a topological
property, but also closely connected to how to obtain the
fundamental cycle. Something we will use since we also have the
following result by Artin \cite{artin}: 
\begin{theorem}\label{zminequqlzmax}
If $(X,0)$ is a rational singularity, then $Z_{max}=Z_{min}$.
\end{theorem} 

This result together with Porposition \ref{thickthinclassyficationofMC}
 shoows the equivalence 
between Theorem \ref{maintheorem} and Theorem \ref{nonmetricversion}. 

Another important property of rational singularities due to Tjurina
\cite{tjurina} is that they are \emph{absolutely isolated}, that is
they can be
resolved by just using blow-up of points. In particular the minimal
good resolution factors through the blow-up of the origin, and hence
resolves the base points of the generic hyperplane sections.

A corollary to this, Theorem \ref{zminequqlzmax} and Proposition
\ref{thickthinclassyficationofMC} is that for a rational surface
singularity being metrically conical is a topological property.
\begin{corollary}\label{topologocaldetermined}
Let $(X,0)$ be a rational singularity that is metrically conical. Then
any normal surface singularity homeomorphic to $(X,0)$ is metrically conical.
\end{corollary}

\subsection{Weighted Homogeneous Surface Singularities.} A complex
surface singularity $(X,0)$ is weighted homogeneous if there is a good
effective $\C^*$ action on $(X,0)$ with $0$ as its only fixpoint. If
$(X,0)\subset (\C^n,0)$ then this is the same as there exist positive
integers $(w_1,\dots,w_n)$ and an action $t\cdot
(z_1,\dots,z_n)=(t^{w_1}z_1,\dots t^{w_n}z_n)$ on $(\C^n,0)$ such that
$f(t\cdot x)=t^\ell f(x)$ for any holomorphic function vanishing on
$(X,0)$. 

Notice that the action of $\C^*$ on $(X,0)$ induces an effective action of $S^1$
on the link of $(X,0)$ and hence the link is Seifert fibred, see for
example Theorem 2.1 of \cite{brandies}. The other direction is
partially true. That is given a star-shaped plumbing graph $\Gamma$, then there
exists a weighted homogeneous singularity with $\Gamma$ as its dual
resolution graph. This follows from the analytic classification of
weighted homogeneous surface singularities
\cite{ConnerRaymonholomorphicseifertfibrations,OrlikWagreichalgbraicsurfacswithkstarcation,pinkhamcstaraction,geometryss}
see also Example 5.1.31 in \cite{nemethinormalsurfacesing}.

By Corollary \ref{topologocaldetermined} we can assume that $(X,0)$ is
weighted homogeneous, but be aware it could be that many of the
analytic types supported on the same toploligical type are not
weigthed homogeneous. This assumption will give us
several advantages.

Assume that the Seifert Invariants of the link of
$(X,0)$ are $\bigl(0,(1,-a),(p_1,q_1),\dots,(p_n,q_n)\bigr)$. Pinkham
introduced the following function in \cite{pinkhamcstaraction}:
\begin{align}
N(\ell) := \ell a -\sum_{i=1}^n\lceil \ell q_i/p_i\rceil \label{pinlhamsfunction}
\end{align} where $\lceil r \rceil$ is the smallest integer larger
than or equal to $r$ and $\ell$ is a non negative integer. We will
call this function the \emph{Pinkham function}.

He then shows that 
\begin{proposition}[corollary 5,8 \cite{pinkhamcstaraction}]
$(X,0)$ is rational if and only if, the graph is
a line or $N(\ell) > -2$ for all $\ell$.
\end{proposition}

It follows by comparing Pinkham's and Laufer's construction, which both
rely on Artin's results, that $N(\ell) = -Z_\ell\cdot E$. So we can
conclude that if $(X,0)$ is rational then $N(\ell) = -1$ for $\ell >0$
and $Z_\ell < Z_{min}$. This is a criterion we will use extensively.

\section{Sufficiency}\label{suffeciency}

We will in the section prove that a singularity with the Seifert
invariants given in Theorem \ref{maintheorem} is rational and
metrically conical. We will first do the case \eqref{singlecurvecase}
with a simple argument and then for \eqref{threeormoreseifertfibrescase}
do this by calculating $N(\ell)$
introduced in \eqref{pinlhamsfunction}. We remarked at the end of last
section how to use this to see if $(X,0)$ is rational, so we need some
result to see how this fact will help us with implying metrically
conicalness of $(X,0)$.

\begin{proposition}
If the exceptional divisor of the minimal good resolution consists of a
single rational curve, then $(X.0)$ is rational and metrically conical.
\end{proposition}

\begin{proof}
It follows from Corollary \ref{linegraph} that $Z_1=Z_{min}$ and therefore by
Laufer's criterion \ref{laufercriterion} that $(X,0)$ is rational. Since the
singularity is rational and hence absolutely isolated by
\cite{tjurina}, the minimal good resolution resolves the
base point of the generic hyperplane sections. So there is only
one vertex and hence one $\mathcal{L}$-node in the resolution graph
determining the thick-thin decomposition. Therefore, $(X,0)$ is metrically
conical by Proposition \ref{thickthinclassyficationofMC}.
\end{proof}

To calculate $N(\ell)$ in case  \eqref{threeormoreseifertfibrescase} we
will study the resolution diagram closer. By a
\emph{bamboo} we mean a valence one vertex together with a string of
valence two vertices. So the resolution graph consists of the central
vertex and $n$ bamboos. Consider a maximal bamboo of the form  
$$\splicediag{8}{30}{
  &&&&\\
 \overtag\Circ{-e_1}{8pt}\lineto[r] &\overtag\Circ{-e_2}{8pt}\dashto[r] &
\overtag\Circ{-e_{k-1}}{8pt} \lineto[r] &\overtag\Circ{-e_k}{8pt}
\lineto[r]& \\ &&&&}$$where the right most edge connects to the central vertex.
To the bamboo we associate a continued fraction $p/q =
[e_k;e_{k-1};\dots;e_2;e_1]$, notice that $(p,q)$ is the Seifert
invariant associated to the singular fibre corresponding to the bamboo.
To all the sub bamboos we associate continued fractions $p_i/q_i =
[e_i;e_{i-1};\dots; e_2;e_1]$. Notice that $p=p_k$, $q=q_k$,
$q_i=p_{i-1}$ for $i>1$ and $q_1=1$.

\begin{lemma}\label{multeq}
Let $Z=mE+\sum_{i=1}^km_{i}E_{i} + Z' \in\z_{top}$, where $Z'$ only
depends on vertices on other bamboos. Then $Z\cdot E_{i}=0$ for all
$i=1,\dots,k$ if and only if $m=m_{1}p$ and $m_{i}=m_{1}p_{i-1}$ for all $i$.
\end{lemma}

\begin{proof}
$Z\cdot E_{i} = 0$ for all $i$ if and only if the multiplicities
satisfy the following set of equations:
\begin{align}
m_{i}e_{i}&=m_{i-1}+m_{i+1}\ \text{ for } 1<i<k-1\label{eq1}\\
m_{1}e_{1}&=m_{2}\label{eq2}\\
m_{k}e_{k}&=m+m_{k-1}.
\end{align}

We want to see that  $m_{i}=m_{1}p_{i-1}$ if and only if these
equations are satisfied. We do this by induction on $i$ starting with
$i = 2$. Using \eqref{eq2}, we get   
\begin{align}
m_{2}&=m_{1}e_{1}=m_{1}[e_{1}]=m_{1}p_{1}.
\end{align}
Assume that $m_{j}=m_{1}p_{j-1}$ for $j<i$. Then from \eqref{eq1} we get that
\begin{align}
m_{i}&=m_{i-1}e_{i-1}-m_{i-2}=m_{1}p_{i-2}e_{i-1}-m_{1}p_{i-3}.\label{eq3}
\end{align}
But 
\begin{align}
[e_{i-1};e_{i-2};\dots;e_{1}] &=
e_{i-1}-\frac{1}{[e_{i-2};\dots;e_{1}]} =
e_{i-1}-\frac{q_{i-2}}{p_{i-2}},
\end{align}
and $q_{i-2}=p_{i-3}$, hence
\begin{align}
[e_{i-1};e_{i-2};\dots;e_{1}]&=\frac{p_{i-2}e_{i-1}-p_{i-3}}{p_{i-2}}.\label{eq4}
\end{align}
Combining equation \eqref{eq3} and \eqref{eq4} gives that
$m_{i}=m_{1}p_{i-1}$. If we consider $m$ as $m_{k+1}$ we also get that
$m=m_{1}p$. For the other direction \eqref{eq2} is trivial and the others
follow from multiplying \eqref{eq4} by $m_{1}p_{i-2}=m_1q_{i-1}$.
\end{proof}

If we have a bamboo with continued fraction $p_i/(p_i-1)$, which is the same
as all the self intersection numbers being $-2$. Then if we run
Laufer's algorithm and get $Z_{p_i}=Z_{min}$, it implies that we have
increased the multiplicity of $E_{ik_i}$ at each step except the
last, but if Laufer's algorithm stops before $Z_{p_i}$ then the above
Lemma cannot be satisfied. For $(X,0)$ to be rational and metrically
conical we need $Z_{min}\geq Z_p$. But as remarked in the end of last
section, this is the same as $N(\ell) = -1$ for $0<\ell<{p_i}$.

\begin{proposition}
Suppose that $(X,0)$ is a weighted homogeneous singularity with
Seifert invariants
$\bigl(0,(1,-(n-1)),(p,p-1),\dots,(p,p-1),(p,p-q),(p,q)\bigr)$
where there are $n\geq 3$ singular fibres. Then $N(\ell) = -1$ for $0<\ell <p$.
\end{proposition}

\begin{proof}
This is a simple calculation:
\begin{align*}
N(\ell) &= \ell (n-1)-\lceil \ell q/p\rceil - \lceil \ell (p-q)/p\rceil -
\sum_{i=1}^{n-2}\lceil \ell (p-1)/p\rceil \\ 
&= \ell (n-1)-\lceil \ell q/p\rceil - \ell - \lceil -\ell q/p\rceil -
\sum_{i=1}^{n-2}\bigl(\ell + \lceil -\ell /p\rceil \bigr) \\ &=
\ell(n-1) -\ell -\bigl(\lceil \ell q/p\rceil + \lceil -\ell
q/p\rceil \bigr) -\ell(n-2)\\
&= -1
\end{align*} where the last equality follows from $ \lceil -\ell
q/p\rceil = -\lfloor \ell q/p\rfloor$ and $ \lceil \ell
q/p\rceil =  \lfloor \ell q/p\rfloor + 1$ since $\gcd(p,q)=1$ and hence
$\ell q/p$ is never an integer.
\end{proof}

\begin{corollary}
Suppose that $(X,0)$ is a weighted homogeneous singularity with
Seifert invariants
$\bigl(0,(1,-(n-1)),(p,p-1),\dots,(p,p-1),(p,p-q),(p,q)\bigr)$
where there are $n$ singular fibres. Then $(X,0)$ is rational and
metrically conical.
\end{corollary}

\begin{proof}
The proposition above takes care of all the cases for $\ell < p$ which is all
that is needed for the metrically conicalness, provided
$Z_{min}=Z_{max}$. For $\ell \geq p$ we write $\ell = \ell' +kp$ with
$l'<p$ and get
\begin{align*}
N(\ell) &= \ell' (n-1) +kp(n-1) 
-\sum_{i=1}^{n}\bigl(\lceil \ell' q_i/p\rceil+ kq_i\bigr) \\ 
&= N(\ell') +kp(n-1) -k(p-1)(n-2) -kq -k(p-q)\\ 
&= N(\ell') + k(n-2).
\end{align*}But $n\geq 3$, $k\geq 1$ and $N(\ell') = 0$ or $-1$. So
$N(\ell) \geq 0$ and hence $(X,0)$ is rational. Notice also that
$N(p)= n-2$ and hence $Z_{min} = Z_p$.
\end{proof}

\section{necessity}\label{necessity}
We will from now on assume that $(X,0)$ is rational and metrically
conical. Let us first consider what happens if there are less than $3$
singular fibres in the Seifert fibration.

\begin{proposition}
Let $(X,0)$ be a rational and metrically conical complex surface
singularity with less than $3$ singular fibres in its Seifert
fibration. Then the exceptional divisor of the minimal good resolution
graph consists of a single smooth rational curve.
\end{proposition}

\begin{proof}
The resolution graph of the minimal good resolution is a line since
there are less than $3$ singular fibres. Then by Corollary \ref{linegraph}
$Z_{min} =Z_1$ and hence $Z_{min}\cdot E < 0$ for any curve $E$
  corresponding to a vertex of valence $1$, because $E^2\leq -2$. This
  implies that if there are more than one curve, and hence two
  vertices of valence one, then there are more than one
  $\mathcal{L}$-node. So $(X,0)$ cannot be metrically conical by
  proposition \ref{thickthinclassyficationofMC}.
\end{proof}

So from now on we also assume that there are at least $3$
singular fibres in the Seifert fibration of the link of $(X,0)$. 
Then following the last section we can assume that $Z_{min}$
satisfies Lemma \ref{multeq} and that $N(\ell) = -1$ for all $\ell >
0$ and $Z_\ell < Z_{min}$. We will need to extract the conditions on
the Seifert invariants from this. Let us first determine the self
intersection number of the central curve. 

\begin{lemma}
Let $(X,0)$ be a metrically conical rational singularity with more than
one curve in its good resolution graph. Assume that $E$ is the central
curve and $a=-E^2$, then $a=n-1$.
\end{lemma}

\begin{proof}
First notice that by running Laufer's Algorithm \ref{laufersalgorithm}
we have that $Z_1 = E + \sum E_{ij}$, and hence 
\begin{align*}
Z_1\cdot E = -a +\sum E_{ij}\cdot E = -a + n.
\end{align*}If $a <n - 1$ then $Z_1\cdot E < -1$ and $(X,0)$ is not
rational by \ref{laufercriterion}. If on the other hand $a\geq n$
then Laufer's 
Algorithm terminates in the first step, since it is clear that for any
valence one and two vertices that $E_{ij}\cdot Z_1 \leq 0$. But this
implies that for any valence one vertex the intersection is strictly
negative, since their weights are less than $-1$. This also implies
that the singularity is rational by Laufer's Criterion.

So the strict transform of $Z_{max}$ intersects the curves corresponding to the
valence one vertices. If there is more than one vertex this contradicts
that they only intersect the central node for a metrically conical
singularity. Hence $a=n-1$.
\end{proof}

\begin{proposition}\label{bamboomultincrease}
Let $(X,0)$ be a rational singularity with Seifert invariants
$\bigl(0,(1,-(n-1)),(p_1,q_1),\dots,(p_n,q_n)\bigr)$. When going
from $Z_k$ to $Z_{k+1}$ then the multiplicity of all exceptional
curves adjacent to the central curve except one is increased by $1$
unless $Z_{k+1} = Z_{min}$ or $k=0$. 
\end{proposition}

\begin{proof}
It is clear that when going from $Z_k$ to $Z_{k+1}$ then each $E_{i1}$
can be added at most one time, and hence 
\begin{align*}
-1 &= Z_{k+1}\cdot E = \bigl(Z_k + E +\sum E_{i1}\bigr)\cdot E
= Z_k\cdot E +E^2 + \sum E_{i1}\cdot E = -1 -(n-1)+ \sum 1
\end{align*}where the sum is taken over all the $i$'s such that
$E_{i1}$ is added.
\end{proof}

\begin{corollary}\label{qigreaterthanphalves}
If $(X,0)$ is a rational singularity where $Z_{min} > Z_2$ with
Seifert invariants
$\bigl(0,(1,-(n-1)),(p_1,q_1),\dots,(p_n,q_n)\bigr)$. Then $q_i \geq
p_i/2$ for all $i$ except one where $q_j< p_j/2$. Equality only
occurs if $p_i=2$.
\end{corollary}

\begin{proof}
Since $Z_2$ is not $Z_{min}$ then by Proposition
\ref{bamboomultincrease} it follows that all but one of the
multiplicities on the $E_{i1}$'s increase when going from $Z_1$ to
$Z_2$. But to increase at that step we have
\begin{align*}
0< (Z_1+E)\cdot E_{ij} =
(2E+\sum E_ij)\cdot E_{i1} = 2E\cdot E_{i1} + E^2_{i1} + E_{i2}\cdot
E_{i1} =  E^2_{i1} + 3.
\end{align*}This implies that $ E^2_{i1} = -2$ and therefore that
$p_i/q_i \leq 2$ and equality means that $q_i$ divides $p$ and hence
$q_i=1$.

The last curve must by the same calculation have $E_{i1}^2 < -2$ and
therefore satisfy the other inequality.
\end{proof}

We will continue with some consequences of Lemma \ref{multeq}.

\begin{proposition}\label{multinmincase}
If $Z_{max}=Z_{min}$ and $(X,p)$ is metrically conical, then $m=\lcm_i(p_i)$.
\end{proposition}

\begin{proof}
If $(X,p)$ is metrically conical, then $Z_{max}\cdot E_{ij}=0$ for all
$i,j$. Since from Lemma \ref{multeq} $m=m_{ik_i}p_i$ for all $i$
(notice that we have turned the numbering of the vertices on the
bamboos around compared to the Lemma), then
there exists a positive integer $l$ such that $m=l\lcm p_i$.  Assume
that $l\neq 
1$. It follows that $l\mid m_{ik_i}$ for all $i$ since
$m_{ik_i}=m/p_i=l\lcm p_i/p_i$. It also follows since $m_{ik_i}\mid m_{ij}$ for
all $i,j$ that $l\mid m_{ij}$ for all $i,j$ by Lemma \ref{multeq}. We can then
define a divisor $Z=(m/l)E+\sum_{i=1}^n\sum_{j=1}^{k_i}(m_{ij}/l)
E_{ij}$ with positive integer coefficients. Since
$Z=\tfrac{1}{l}Z_{min}$ it implies $Z\in\z_{top}$. But $Z<Z_{min}$,
hence we get a contradiction.
\end{proof}

\begin{lemma}\label{multincrease}
If the continued fraction associated to a bamboo is $p/q$, then we do
not increase the multiplicity on vertices along the bamboo when we go
from $Z_{lp-1}$ to $Z_{lp}$.
\end{lemma} 

\begin{proof}
By our definition of $Z_k$ we have that $Z_{lp-1}\cdot E_i \leq 0$ for
all the exceptional curves along the bamboo. By Lemma \ref{multeq} (or
more precisely its proof) we have that $Z_{lp}\cdot E_i = 0$ for all
$i$. For any multiplicity along the bamboo to be increased when going
from $Z_{lp-1}$ to $Z_{lp}$, we would need to increase the multiplicity of
$E_1$ (the closest to the central vertex). So $Z_p = Z_{p-1} + E + E_1
+ aE_2+...$ where $a=0$ or $1$ depending on wether $E_2$ needs to be added.
Hence 
\begin{align*}
0=Z_{lp}\cdot E_1 = Z_{lp-1}\cdot E_1 + E\cdot E_1 + E_1^2 + aE_2
\cdot E_1 =  Z_{lp-1}\cdot E_1 +1 +E_1^2 +a 
\end{align*}If $a=0$ we get a contradiction since $ Z_{lp-1}\cdot E_1
\leq 0$ and $E_1\leq -2$. If $a=1$ then one can do the same for
$E_p\cdot E_2$, either reaching a contradiction or concluding that
$E_3$ has been added. Continuing this process we either get a contradiction
or we end up adding all the curves along the bamboo. But then doing
the calculation for the last curve gives
\begin{align*}
0=Z_{lp}\cdot E_k = Z_{lp-1}\cdot E_1 + E_{k-1}\cdot E_k + E_k^2
\cdot E_1 =  Z_{lp-1}\cdot E_1 +1 +E_1^2 
\end{align*}and we get a contradiction from $ Z_{lp-1}\cdot E_k
\leq 0$ and $E_k\leq -2$.
\end{proof}

\begin{corollary}\label{difrenetbamboos}
If two bamboos have the associated continued fractions $p/q_1$ and
$p'/q_2$ with $p'$ a divisor of $p$, then $Z_{min} \leq Z_p$. 
\end{corollary}

\begin{proof}
By the Lemma none of the multiplicities along those
bamboos are increased when going from $Z_{p-1}$ to $Z_p$, which by Proposition
\ref{bamboomultincrease} implies that Laufer's Algorithm stops at
$Z_p$ if it had not stopped before.
\end{proof}

\begin{corollary}\label{seifeertformbeforeequilityofdegrees}
If $(X,0)$ is rational and metrically conical with Seifert invariants
$\bigl(0,(1,-(n-1)),(p_1,q_1),\dots,(p_n,q_n)\bigr)$ with $n\geq
3$. Then $m=\lcm(p_i,p_j)$ for all $i\neq j$. 

Moreover we can organize
the pairs so that we have
$m=\rho\prod_{i=1}^{n'} \rho_i$, $p_j =\rho\prod_{i=1,\ i\neq j}^{n'} \rho_i$
for $j\leq n'$, $p_j = m$ for $j> n'$ and $m_i = \rho_i$ where
the $\rho_i$'s are pairwise coprime integers.
\end{corollary}

\begin{proof}
Since $m$ is the multiplicity of the central node in $Z_{min}$ and
$(X,0)$ is metrically conical it follows from Proposition
\ref{multinmincase} that $Z_{min}=\lcm(p_1,\dots,p_n) = m$ for all
$i$. Now assume that there are $p_i$ and $p_j$ such that
$p' := \lcm(p_i,p_j) < m$. Then $p'=l_ip_i = l_jp_j$ and by Lemma
  \ref{multincrease} the multiplicity does not increase along both
  bamboos when we go from $p'-1$ to $p'$. This implies that Laufer's
  algorithm would have to stop at $p'$ if not before, which contradicts that
  $p'<m$.

By definition of $m_i$ we have $m=m_ip_i$. Now assume that for some
$m_j$ we have $m_{ij} :=\gcd(m_i,m_j) > 1$, then $(m_i/m_{ij}) p_i =
(m_j/m_{ij})p_j$ is a common multiple of $p_i$ and $p_j$ that is less
than $m$. Hence $\gcd(m_i,m_j) = 1$ for all $i\neq j$, and we can
write
\begin{align*}
m=m_1\cdots m_n\rho
\end{align*}for some $\rho$. Ordering the $p_i$'s so that the ones
with $m_i=1$ are the last $n-n'$ we get the factorization of $m$ and
$p_i$ by setting $\rho_i=m_i$.
\end{proof}

\begin{proposition}\label{sumofqis}
Let $(X,0)$ be a rational metrically conical singularity with at least
$3$ singular fibres, Seifert
invariants $\bigl(0,(1,-(n-1)),(p_1,q_1),\dots,(p_n,q_n)\bigr)$ and
let $m$ be the multiplicity of the central curve in $Z_{min}$. Then
the following equation is satisfied:
\begin{align}
 (m-1)(n-1)+1 = \sum_{i=1}^nm_iq_i\label{multequation}
\end{align} where $m_i=m/p_i$.
\end{proposition}

\begin{proof}
Since $Z_{min}=Z_m$ that means we have to increase the multiplicity of
the central curve $m$ times. The first time we also increase the
multiplicity of all the neighbouring vertices. The next $m-2$ times we
increase the multiplicity of all but one of the neighbouring vertices, and
the last increase does not increase any neighbouring vertices. This
follows from Proposition \ref{bamboomultincrease} and
\ref{multincrease}. Summing the numbers of increments gives the left
hand side of the equation.

The multiplicity of $E_{i1}$ is $m_iq_i$ by Lemma \ref{multeq}, and
hence has been increased that many times. Summing these numbers over
all the bamboos gives us the right hand side.

Since both these count the number of times an $E_{i1}$ has been
increased from going from $Z_0$ to $Z_{min}$ the sides have to be
equal.
\end{proof}
An alternative proof not using the combinatorial logic behind the
plumbing graph can be given as follows. Use the form of the Seifert
invariants obtained in Corollary
\ref{seifeertformbeforeequilityofdegrees}, then use that $N(\ell) = -1$
for all $\ell\in 1,\dots,p-1$. If one calculates $N(p-1)$ and puts it
equal to $-1$, then one will get the same formula.

\begin{corollary}\label{commonseifertdegree}
If $(X,0)$ is rational and metrically conical with Seifert invariants
$\bigl(0,(1,-(n-1)),(p_1,q_1),\dots,(p_n,q_n)\bigr)$ with $n\geq
3$. Then $p_i=p_j=m$ for all $i,j$. 
\end{corollary}

\begin{proof}
We will use the form of the Seifert invariants obtained in Corollary
\ref{seifeertformbeforeequilityofdegrees}, and the fact that $N(\ell)=-1$
for all $\ell\in 1,\dots,p-1$. Let us first calculate $N(\ell)+N(m-\ell)$ for
arbitrary $\ell$.
\begin{align}
N(\ell)+N(m-\ell) &= (n-1)\ell-\sum_{i=1}^n \Big\lceil \ell\frac{q_i}{p_i} \Big\rceil +
(n-1)(m-\ell)-\sum_{i=1}^n \Big\lceil (p-\ell)\frac{q_i}{p_i}
\Big\rceil  \label{sumofdualN(k)} \\
 &= (n-1)m-\sum_{i=1}^n \Big\lceil \ell\frac{q_i}{p_i} \Big\rceil 
-\sum_{i=1}^n \Big\lceil \rho_iq_i -\ell\frac{q_i}{p_i} \Big\rceil
\nonumber \\
 &= (n-1)m-\sum_{i=1}^n \biggl(\Big\lceil \ell\frac{q_i}{p_i}
 \Big\rceil + \Big\lceil -\ell\frac{q_i}{p_i} \Big\rceil \biggr)
-\sum_{i=1}^n \rho_iq_i  \nonumber \\
 &= (n-1)m-\sum_{i=1}^n \biggl(\Big\lceil \ell\frac{q_i}{p_i}
 \Big\rceil - \Big\lfloor \ell\frac{q_i}{p_i} \Big\rfloor \biggr)
-\big((n-1)(m-1) +1\big)  \nonumber \\
 &= n-2-\sum_{i=1}^n \biggl(\Big\lceil \ell\frac{q_i}{p_i}
 \Big\rceil - \Big\lfloor -\ell\frac{q_i}{p_i} \Big\rfloor \biggr),
\nonumber
\end{align}where we use Proposition \ref{sumofqis} for the second to
last equality.

Now the difference between the ceiling and the floor is
\begin{align*}
\Big\lceil \ell\frac{q_i}{p_i}
 \Big\rceil - \Big\lfloor -\ell\frac{q_i}{p_i} \Big\rfloor
 = \begin{cases}
1 &\text{if }\ell\frac{q_i}{p_i}\text{ is not an integer}\\
0 &\text{if }\ell\frac{q_i}{p_i}\in \Z.
\end{cases}
\end{align*}Since $\gcd(p_i,q_i) = $ we have that $\ell\frac{q_i}{p_i}$
is only an integer if $p_i$ divides $\ell$. So if we define
\begin{align*}
D(\ell) = \# \{p_i\ \vert\ p_i\text{ does not divide }\ell\},
\end{align*}and use that $N(l) = -1$ for all $\ell\in 1,\dots m-1$ if
$(X,0)$ is rational and metrically conical, \eqref{sumofdualN(k)} turns into
\begin{align*}
D(\ell) = n
\end{align*}for all $\ell$. Consequently $p_i$ does not divide $\ell$ for any
value of $\ell$.

Assume that there exists $\rho_i> 1$, then $p_i = m/\rho_i < m$ and
hence we can choose $\ell=p_i$ which clearly implies that $p_i$ divides
$\ell$. Hence $\rho_i=1$ for all $i$ and therefore $p_i=m/\rho_i = m$ for
all $i$.
\end{proof}

Combining what we have so far gives us the following.
\begin{proposition}\label{summerary}
Let $(X,0)$ be a metrically conical rational surface singularity, and
assume that it has $n>2$ singular fibres. Then the Seifert invariants
are
$\bigl(0,(1,-(n-1)),(p,q_1),\dots,(p,q_{n-1}),(p,q_n)\bigr)$. Moreover
$q_i> p/2$ except for one where the other
inequality is satisfied and $(p-1)(n-1)+1 =\sum_{i=1}^{n}q_i$. 
\end{proposition}

\begin{proof}
This is simply combining Corollary \ref{qigreaterthanphalves},
proposition \ref{commonseifertdegree} and equation \eqref{multequation}.
\end{proof}

We will now use this to calculate the Pinkham function for all
$0<\ell<p$ and reach some conclusions on the $q_i$'s. 

\begin{lemma}
Assume that \eqref{multequation} is satisfied. Then $N(\ell)=-1$ is
equivalent to
\begin{align}
\bigg\lceil \ell\frac{n-2+q_1+\cdots+q_{n-1}}{p}\bigg\rceil
= \sum_{j=1}^{n-1} \lceil \ell q_j/p \rceil.\label{ceilingequation}
\end{align} for all $\ell < p$.
\end{lemma}

\begin{proof}
We will use \eqref{multequation} to replace $q_n$ in the calculation
of $N(\ell)$. So
\begin{align*}
 N(\ell) &=\ell (n-1) -\sum_{j=1}^{n}
 \lceil \ell q_j/p \rceil =\ell (n-1) -\lceil \ell q_n/p \rceil
 -\sum_{j=1}^{n-1} \lceil \ell q_j/p \rceil \\
&= \ell (n-1) - \bigg\lceil
\ell\frac{p(n-1)-n+2-q_1-\cdots-q_{n-1}}{p}\bigg\rceil  -\sum_{j=1}^{n-1}
\lceil \ell q_j/p \rceil \\
&= \ell (n-1) - \bigg\lceil
\ell(n-1) -\ell\frac{n-2+q_1+\cdots+q_{n-1}}{p}\bigg\rceil  -\sum_{j=1}^{n-1}
\lceil \ell q_j/p \rceil \\
&= \ell (n-1) -
\ell(n-1)  -
\bigg\lceil-\ell\frac{n-2+q_1+\cdots+q_{n-1}}{p}\bigg\rceil
-\sum_{j=1}^{n-1} \lceil \ell q_j/p \rceil \\
&=  -
\bigg\lceil-\ell\frac{n-2+q_1+\cdots+q_{n-1}}{}\bigg\rceil
-\sum_{j=1}^{n-1} \lceil \ell q_j/p \rceil
\end{align*} Since the ceiling of a negative number is the same as minus
the floor of that number we get that the $N(\ell)=-1$ is equivalent to
\begin{align*} 
\bigg\lfloor \ell\frac{n-2+q_1+\cdots+q_{n-1}}{p}\bigg\rfloor+1
= \sum_{j=1}^{n-1} \lceil \ell q_j/p \rceil 
\end{align*}

Now $\gcd(p,q_n)=1$ and hence $\ell q_n/p$ 
is not an integer as long as $\ell < p$. This means that what is
inside the first floor function is
never an integer because we just remove an integer from $\ell
q_n/p_n$. This then implies that the ceiling is the floor $+1$ and
equation \eqref{ceilingequation} follows.
\end{proof}

Let us introduce some notation to make things easier to read. First
set  $a_i := \lceil \ell q_i/p \rceil$ for $i=1,\dots, n-1$, and
$b=\sum_{i=1}^{n-1}a_i$. Let $\{q/p\}=q/p-\lfloor q/p \rfloor$ be the
fractional part of a rational number.

Set $\epsilon_i = p-q_i$ for $1<i< n$ and $\epsilon =
\epsilon_2+\dots \epsilon_{n-1}$. Notice that by our assumptions on
the $q_i$ we have that $\epsilon_i < p/2$ for all $i$ except maybe one.

\begin{lemma}\label{fractionalinequalities}
If $(X,0)$ is a rational and metrically conical singularity with more
than $2$ singular fibres. Then
\begin{align}
p\{\ell q_1/p\} \in \Bigl(\ell\bigl(\epsilon
-(n-2)\bigr),p\Bigr)\label{groupin}
\end{align}for all $0<\ell < \min\{p/\epsilon_2,\dots,p/\epsilon_{n-1}\}$.
\end{lemma}

\begin{proof}
Using this notation introduced before the Lemma we see that
\begin{align*}
\bigg\lceil \ell\frac{n-2+q_1+\cdots+q_{n-1}}{p}\bigg\rceil &=
\Big\lfloor \ell \frac{q_1}{p}\Big\rfloor +\dots+ \Big\lfloor \ell
\frac{q_{n-1}}{p}\Big\rfloor
  +\bigg\lceil \Bigl\{\ell \frac{q_1}{p}\Bigr\}+\dots+\Bigl\{\ell
  \frac{q_{n-1}}{p} \Bigr\} 
  +\ell \frac{n-2}{p}
  \bigg\rceil \\
&= b - (n-1) +\bigg\lceil \{\ell q_1/p\}+\dots+\{\ell q_{n-1}/p\} +\ell (n-2)/p
  \bigg\rceil.
\end{align*}

Since $\bigg\lceil \ell\frac{n-2+q_1+\cdots+q_{n-1}}{p}\bigg\rceil =
b$ by \eqref{ceilingequation} it follows that
\begin{align*}
n-2 < \{\ell q_1/p\} +\dots+\{\ell q_{n-1}/p\} +\ell(n-2)/p \leq n-1.
\end{align*}Notice that the first inequality is trivial if $\{\ell
q_1/p\} +\dots+\{\ell q_{n-1}/p\} > n-2$ and the second is trivial if $\{\ell
q_1/p\} +\dots+\{\ell q_{n-1}/p\} \leq 1$.

We have that $\{\ell
q_i/p\} = 1 - \{\ell\epsilon_i/p\}$. Now if $\ell\epsilon_i < p$ then
we can remove the fractional part on the right side of the
equation. So from now on we assume that $\ell < p/\epsilon_i$ for all
$i$. 

We will rewrite the above equation using this notation and hypothesis:
\begin{align*}
n-2 < \{\ell q_1/p\}+1-\ell\epsilon_2/p +\dots+1 -
\ell\epsilon_{n-1}/p  +\ell(n-2)/p \leq n-1,
\end{align*}which is the same as 
\begin{align*}
\ell\epsilon/p - \ell(n-2)/ p< \{\ell q_1/p\} 
\leq \ell\epsilon/p -\ell(n-2)/p +1.
\end{align*}

We also have that
$\epsilon \geq n-2$ since each $\epsilon_i \geq 1$. So the right hand
side can be replaced by $ < 1$ since
$\{\ell q_1/p\}$ is always smaller than $1$. Notice this
implies that 
\begin{align*}
\{\ell q_1/p\} \in \Bigl(\ell\epsilon/p  -\ell(n-2)/p,1\Bigl)
\end{align*}or
\begin{align*}
p\{\ell q_1/p\} \in \Bigl(\ell\bigl(\epsilon  -(n-2)\bigr),p\Bigr).
\end{align*}
\end{proof}

\begin{proposition}\label{q1equalpminus1}
Let $(X,0)$ be a rational metrically conical singularity and assume
that at least one of the $\epsilon_i > 1$ then $q_1=p-1$.
\end{proposition}

\begin{proof}
Order the indices such that
$\epsilon_{n-1} = \max{\epsilon_i}$. Then $\epsilon-(n-2)>
1$ and equation \eqref{groupin}
is satisfied
 for all $\ell < \min(p/\epsilon) \leq
p/\epsilon_{n-1} \leq p/2$. Notice that $p\{\ell q_1/p\}$ is the same
as the smallest positive representative of the class of $\ell q_1$ in
$\Z/p\Z$. Hence if we give $\Z/p\Z$ the order inherited from $\Z$ the
condition \eqref{groupin} is the same as saying that the map multipliction by
$q_1$ takes the first $p/\epsilon_{n-1}-1$ non-zero elements of
$\Z/p\Z$ to the last $p/\epsilon_{n-1}-1$ elements. It now follows
from the following proposition that $q_1=p-1$.
\end{proof}

\begin{proposition}
Let $q\in \Z/p\Z$ be an invertible element, and $0 < k< p/2$. Suppose that
the map $\morf{m_q}{\Z/p\Z}{\Z/p\Z}$ given by multiplication by $q$,
takes the first $k$ non-zero elements to the last $k$ non-zero
elements. Then $q=p-1$.
\end{proposition}

\begin{proof}
Because $q$ is invertible $m_q$ is injective. We also have that
$m_q(1)=q$ and hence $q$ is one of the last $k$ elements. Let 
$a=p-q$ it follows that $a$ is one of the first $k$ elements. If
$l$ is one of the first $k$ non-zero elements, then
$m_q(l)= ql= (m-a)l=ml-al$ is one of the last $k$ elements and hence
$al$ is one of the first $k$ elements. So $a$ is an invertible element
such that $m_a$ takes the first $k$ elements to the first $k$
elements. Assume that $b$ is one of the first $k-1$ elements. Then
$m_a(b+1) = m_a(b) + a$. Since both $m_a(b)$ and $a$ are among the first $k$
elements and $k<p/2$, then $m_a(b+1)$ must be less than $p$. Hence
$m_a(b+1) > m_a(b)$ in the order inherited from $\Z$. This implies
that $m_a$ is a bijective map on a finite ordered set that preserves
the ordering, which means that $m_a$ is the identity. So $a=1$ and
therefore $q=p-1$.
\end{proof}

We can now prove the other direction of Theorem \ref{maintheorem}.

\begin{theorem}
Let $(X,0)$ be a rational metrically conical singularity with $n>2$
singular fibres. Then the Seifert invariants of $(X,0)$ are
$\bigl(0,(1,-(n-1),(p,p-1),\dots,(p-1),(p,p-q),(p,q)$ where $\gcd(p,q)=1$.
\end{theorem}

\begin{proof}
The only things that do not follow from Proposition \ref{summerary}
is the values of the $q_i$'s.

Notice that the choice of $q_1$ in Lemma \ref{fractionalinequalities}  and
Proposition
\ref{q1equalpminus1} is arbitrary among $q_1,\dots,q_{n-1}$. Hence the
conclusion is that if one of the $q_i$ is not equal to $p-1$ the rest
are.

So we can assume that $q_1=\dots = q_{n-2} = p -1$. Then it follows from
equation  \eqref{multequation} that $(p-1)(n-1)+1 = q_n
+q_{n-1} +\sum_{i=1}^{n-2}(p-1)$ implies that $p=q_n
+q_{n-1}$.  
\end{proof}

\bibliography{general}

\end{document}